\documentclass[article,12pt]{article}
\usepackage{setspace}
\usepackage{times}
\usepackage{dsfont}
\usepackage{mathrsfs}
\usepackage[colorlinks]{hyperref}
\usepackage{xcolor,soul}
\usepackage{mathscinet}
\usepackage{latexsym}
\usepackage{amsthm}
\usepackage{amssymb}
\usepackage{amsfonts}
\usepackage{amsmath}
\usepackage{longtable}
\usepackage{graphicx}
\usepackage{multirow}
\usepackage{multicol}
\usepackage{booktabs}
\usepackage{fullpage}
\usepackage{cases}
\usepackage{natbib}
\usepackage{enumerate}

\newtheorem{theorem}{Theorem}[section]
\newtheorem{thm}[theorem]{Theorem}

\newtheorem{lem}[theorem]{Lemma}
\newtheorem{proposition}[theorem]{Proposition}

\newtheorem{corollary}[theorem]{Corollary}

\theoremstyle{definition}

\newtheorem{defn}[theorem]{Definition}

\theoremstyle{remark}

\newtheorem{rem}[theorem]{Remark}
\numberwithin{equation}{section}

 \DeclareMathAlphabet{\mathpzc}{OT1}{pzc}{m}{it}

\newcommand{\E}{\mathbb{E}}

\newcommand{\R}{\mathbb{R}}

\newcommand{\PP}{\mathbb{P}}

\newcommand{\Be}{\begin{equation}}
\newcommand{\Ee}{\end{equation}}
\newcommand{\Bs}{\begin{split}}
\newcommand{\Es}{\end{split}}
\newcommand{\Bes}{\begin{equation*}}
\newcommand{\Ees}{\end{equation*}}
\newcommand{\BT}{\begin{thm}}
\newcommand{\ET}{\end{thm}}
\newcommand{\Bp}{\begin{proof}}
\newcommand{\Ep}{\end{proof}}
\newcommand{\BL}{\begin{lem}}
\newcommand{\EL}{\end{lem}}
\newcommand{\BP}{\begin{proposition}}
\newcommand{\EP}{\end{proposition}}
\newcommand{\BC}{\begin{corollary}}
\newcommand{\EC}{\end{corollary}}
\newcommand{\BR}{\begin{rem}}
\newcommand{\ER}{\end{rem}}
\newcommand{\BD}{\begin{defn}}
\newcommand{\ED}{\end{defn}}
\newcommand{\BI}{\begin{itemize}}
\newcommand{\EI}{\end{itemize}}

\title{Sparse Poisson Regression with Penalized Weighted Score Function\thanks{Jinzhu Jia and Fang Xie contributed equally, they are co-first authors and are listed in alphabet orders, Lihu Xu is the corresponding author.}}
\author{Jinzhu Jia\thanks{School of Mathematical Sciences and Center for Statical Science, Peking University, Beijing, China} \and Fang Xie\thanks{1. Department of Mathematics, Faculty of Science and Technology,
University of Macau, Av. Padre Tom\'{a}s Pereira, Taipa Macau, China. 2. UMacau Zhuhai Research Institute, Zhuhai, China. FX is supported by the grants Macao S.A.R FDCT  030/2016/A1, 049/2014/A1 and  NNSFC 11571390 and University of Macau MYRG2015-00021-FST, 2016-00025-FST. LX is supported by the grants Macao S.A.R FDCT  030/2016/A1, 049/2014/A1 and  NNSFC 11571390 and University of Macau MYRG2015-00021-FST, 2016-00025-FST.} \and Lihu Xu\footnotemark[3]}

\begin{document}
\date{}
\maketitle
\begin{abstract}
\indent   We proposed a new penalized method in this paper to solve sparse Poisson Regression problems. Being different from $\ell_1$ penalized log-likelihood estimation, our new method can be viewed as a penalized weighted score function method. We show that under mild conditions, our estimator is $\ell_1$ consistent and the tuning parameter can be pre-specified, which owns the same good property of the square-root Lasso. The simulations show that our proposed method is much more robust than traditional sparse Poisson models using $\ell_1$ penalized log-likelihood method.
\end{abstract}
\ \ \ \ \ \ \ \ \ \

\section{Introduction}

Poisson regression is a special generalized linear model \citep{nelder1972generalized} which is widely used to model count data. Let $(x_{1},y_{1}),\cdots,(x_{n},y_{n})$ be independent pairs of observed data which are realizations of random vector $(X,Y)$,
with $p$-dimensional covariates $X\in \R^p$ and univariate response variable $Y\in \R$. $(X,Y)$ is assumed to satisfy the conditional distribution $Y|X=x\sim {\rm Poisson}(\mu(x))$ with $\log(\mu(x))=x^{T}\beta^{*}$, where $\beta^{*}\in\mathbb{R}^{p}$ is an unknown parameter vector to be estimated.

In this paper, we are concerned with a sparse Poisson regression problem when the number of covariates (or predictors) is much larger than the number of observations, i.e. $p\gg n $,  which is a variable selection (or model selection) problem for high-dimensional data. For linear models, now researchers have developed several methods such as Lasso \citep{tibshirani1996regression}, adaptive Lasso \citep{zou2006adaptive}, SCAD \citep{fan2001variable}  and so on.  Lasso is a very popular method not only due to its interpretability and prediction performance \citep{tibshirani1996regression},  but also because it is a convex problem and can be computed easily and fast \citep{Friedman2010Regularization}. It is well known that when incoherent condition (or irrepresentable condition) holds,  Lasso estimator is sign consistent \citep{zhao2006model,zou2006adaptive,jia2013lasso}. When a restricted eigenvalue condition holds, Lasso estimator can be $\ell_2$ consistent \citep{bickel2008simultaneous}.  Because incoherent condition might not hold, more steps are used to relax this condition such as adaptive Lasso \citep{zou2006adaptive} and Puffer transformation \citep{jia2015preconditioning}. Other non-convex penalized methods like SCAD \citep{fan2001variable} and  MCP \citep{zhang2010nearly} can also be used to study sparse models.

Variable selection problems for generalized linear models have also gained great attentions in recent years. For instance, \cite{ravikumar2010high-dimensional} studied $\ell_1$ regularized logistic regression models, while \cite{Li2015Consistency} analyzed the consistency of $\ell_1$ penalized Poisson regression models. Moreover, \cite{raginsky2010compressed} studied the performance bounds for compressed sensing (CS) under Poisson models, \cite{ivanoff2016adaptive} also considered a data-driven tuning method for sparse and structure sparse functional Poisson regression models.

 It is now well known that  the Lasso problem in linear regression models, a good tuning parameter choice depends on the unknown parameter $\sigma^2$ which is the homogeneous noise variance in linear models \citep{bickel2008simultaneous}. To solve this problem, \cite{belloni2011square-root} proposed square-root Lasso, which alternatively replaces
the original score function in \citep{bickel2008simultaneous} by the square root of this function.
In previous studies \citep{raginsky2010compressed,Li2015Consistency,ivanoff2016adaptive}, variable selections for sparse Poisson models are obtained via penalized loglikelihood methods, which have the same problem as the Lasso problem. Moreover, Poisson noises are not homogeneous any more, a unique penalty for all of the different coefficient is not a good choice  \citep{ivanoff2016adaptive}. In this paper, we propose a new {\it penalized weighted score function} method to study sparse Poisson regression, and show that it gives consistent estimator of the parameters in sparse Poisson models and provides a direct choice for the tuning parameter. The simulations show that our proposed method is much more robust than traditional sparse Poisson models using $\ell_1$ penalized log-likelihood method.

The rest of the paper is arranged as follows. In Section \ref{sec: penalized weighted score function}  we first review square root Lasso and explain why it could be viewed as a penalized weighted score function. Then we apply this idea to sparse Poisson models and propose our method. Section \ref{sec: theories} provides finite sample and asymptotic bounds for our new estimator. In Section \ref{sec:experiments}, we conduct experiments to show the robustness of our method. Section \ref{sec:proofs} gives the detailed proofs for our theoretical results.

%There are two advantages comparing our method with other method solving Poisson regression.
%\indent 1. The penalty coefficient $\lambda$ does not depend on the noise level $\sigma_i$.\\
%\indent 2. We prove the $l_1$-consistency of $\hat{\beta}$.
\section{Penalized weighted score function}
\label{sec: penalized weighted score function}

We first briefly give a few notations used in this paper.
\subsection{Notations}
Let $[n]=\{1,2,\cdots,n\}$, $[p]=\{1,2,\cdots,p\}$ and $[d]=\{1,2,\cdots,d\}$. For any $d$-dimensional vector $v=(v_1,\cdots,v_d)^T$, denote $\|v\|^q_q=\sum\limits_{i=1}^d |v_i|^q$ for any $q\in(0,+\infty)$ and denote $\|v\|_\infty=\max\limits_{i\in[d]}|v_i|$.

Write $X=(x_1,\cdots,x_n)^T\in\R^{n\times p}$ and $Y=(y_1,\cdots,y_n)^T\in\R^n$. Denote by $T=supp(\beta^{*})=\{j\in[p]: \beta^*_j\neq 0\}$ the non-zero coordinate of $\beta^*$ and let $s=|T|$ be the number of non-zero elements of $\beta^*$. $I_n$ denotes for the $n\times n$ identity matrix.

Denote by $\{a_{n}\}_{i=1}^{n}$ and $\{b_{n}\}_{i=1}^{n}$ two sequences, the notation $b_{n}=O(a_{n})$ means that there exists a constant $C>0$ such that $b_{n}\le Ca_{n}$ for all $n\ge 1$ and the notation $b_{n}=o(a_{n})$ means that $\lim\limits_{n\rightarrow \infty}\frac{b_{n}}{a_{n}}=0$.

If $f$ is a function, we denote by $\nabla f$ the gradient of $f$.

\subsection{Square root Lasso revisited}
We now review square root Lasso and treat it as a penalized weighted score function method. The Lasso is defined as follows.
\begin{equation}
\label{eqn:lasso}
\min_{\beta} \|Y - X\beta\|_2^2 + 2 \lambda \|\beta\|_1,
\end{equation}
where $Y\in \mathbb R^{n\times 1}, X \in \mathbb R^{n\times p}$ and $\beta \in \mathbb R^{p\times 1}$. The solution of the Lasso satisfies KKT conditions defined as follows:

\begin{equation}
\label{eqn:KKT}
X^T(Y-X\beta) = \lambda \vec s,
\end{equation}
where $\vec s_j$ is a subgradient of $|\beta_j|$ which is the sign of $\beta_j$ if $\beta_j\neq 0$ and can be any value belonging to $[0,1]$ when  $\beta_j= 0$.  From Equation \eqref{eqn:KKT}, we see $$|X_j^T(Y-X\beta)| \leq \lambda\qquad\qquad \ {\rm for\ all\ }  j=1,\cdots,p.$$

Note that $X^T(Y - X\beta)$ is the score function for linear model with Gaussian noises, that is  $(X,Y)$ follows $$Y = X\beta^* + \epsilon,$$ with $\epsilon \sim N(0,\sigma^2 I_n)$.  To have a good estimator ($\ell_2$ consistent for example), the choice of $\lambda$ should satisfy the condition $\lambda > c \|X^T(Y-X\beta^*)\|_\infty$ for some positive constant $c$ \citep{bickel2008simultaneous}. The score function evaluated at $\beta^*$ is $X^T(Y-X\beta^*)$ which has a multi-normal distribution with mean $0$ and co-variance matrix $\sigma^2 X^TX$, this suggests that the choice of $\lambda$ also depends on the unknown parameter $\sigma^2$. One way to remove this unknown parameter is to use a weighted (or scaled) score function defined as
\begin{equation*}
%\label{eqn:weightscore}
\frac{X^T(Y-X\beta)}{\|Y-X\beta\|_2},
\end{equation*}
whose distribution at $\beta^*$ does not depend on $\sigma^2$ any more.

Setting the penalized weighted score function to be 0, i.e.
\begin{equation}
\label{eqn:weightscore}
\frac{X^T(Y-X\beta)}{\|Y-X\beta\|_2} - \lambda \vec{s} =0,
\end{equation}
leads to the following optimization problem:
\begin{equation}
\label{eqn: square-root Lasso}
\min_{\beta} \|Y - X\beta\|_2 + \lambda \|\beta\|_1,
\end{equation}
which is the square root Lasso defined in \cite{belloni2011square-root}.

This viewpoint of treating square root Lasso as a  penalized weighted score function method could be applied to other problems especially for heteroscedastic models. In the next subsection, we give details on how to apply this idea to sparse Poisson models.

\subsection{Penalized weighted score function method for sparse Poisson regression}
Let $(x_{1},y_{1}),\cdots,(x_{n},y_{n})$ be independent pairs of observed data which are realizations of random vector $(X,Y)$,
with $p$-dimensional covariates $X\in \R^p$ and univariate response variable $Y\in \R$. $(X,Y)$ is assumed to satisfy the conditional distribution $Y|X=x\sim {\rm Poisson}(\mu(x))$ with $\log(\mu(x))=x^{T}\beta^{*}$, where $\beta^{*}\in\mathbb{R}^{p}$ is an unknown parameter vector to be estimated.
 Denoting $x_i=(x_{i1},\cdots,x_{ip})^T$, without loss of generality, we assume
$$\frac1{n}\sum\limits_{i=1}^n{x^{2}_{ij}}=1,\ \ \ \ \ \ \ \ \ \ \ {\rm for\ all}\ j\in[p].$$

Under the above settings, the negative loglikelihood (up to a scale and a constant shift) is defined as follows:
\begin{equation}
\label{eqn:loglikelihood}
\ell(\beta) = \frac{1}{n}\sum_{i=1}^n(-y_i x_i^T\beta + e^{x_i^T\beta}).
\end{equation}

Sparse Poisson regression model could be gained via $\ell_1$ penalized loglikelihood defined as follows:
\Bes
\min\limits_{\beta\in\R^{p}}\big\{ \frac1{n}\sum\limits_{i=1}^{n}(-y_{i}x_{i}^{T}\beta+e^{x_{i}^{T}\beta})+\bar{\lambda}\|\beta\|_{1}\big\},
\Ees
where the penalty level $\bar{\lambda}>0$.  From the KKT optimality condition, we get $\bar{\lambda}\ge \|\nabla \ell(\beta)\|_\infty$. To get a good estimator, usually we require that the inequality $\bar{\lambda}\ge c \|\nabla \ell(\beta^{*})\|_\infty$  for some constant $c>1$ holds with high probability \citep{Li2015Consistency}. However, for the score function valued at $\beta=\beta^{*}$:
$$\nabla \ell(\beta^{*})=-\frac1{n}\sum\limits_{i=1}^{n}x_{i}(y_{i}-e^{x_{i}^{T}\beta^*}),$$
the random part $y_{i}-e^{x_{i}^{T}\beta^*}$ has variance $e^{x_{i}^{T}\beta^{*}}$, which is also the rate parameter of the Poisson random variable $Y_{i}|X_{i}=x_{i}$. If the rate parameter is very large, the penalty coefficient $\lambda$ will be very sensitive.

In this paper, we apply the idea from square root Lasso to solve the above problem on choosing penalty level, let us briefly introduce
our new method as follows. We give weights to the items in score function and develop an $\ell_1$ penalized weighted score function method, which solves the following equation:
\begin{equation}
\label{weightedL1scoreFunction}
-\frac1{2n}\sum\limits_{i=1}^{n}\frac{x_{i}(y_{i}-e^{x_{i}^{T}\beta})}{\sqrt{e^{x_{i}^{T}\beta}}} + \lambda \vec s = 0,
\end{equation}
where $\vec s$ is the sub-gradient of $\|\beta\|_1$. By a careful derivation, we found that the solution to Equation \eqref{weightedL1scoreFunction} is equivalent to solve the following convex optimization problem:
\begin{equation}
\label{eqn:mainEQ}
\hat\beta = \arg\min_{\beta\in \mathbb R^p}\left\{\frac1{n}\sum\limits_{i=1}^{n}(y_{i}e^{-\frac1{2}x_{i}^{T}\beta}+e^{\frac1{2}x_{i}^{T}\beta}) + {\lambda}\|\beta\|_1 \right\},
\end{equation}
where $\lambda>0$ is a new penalty level which will not depend on any rate parameter $e^{x_{i}^{T}\beta^{*}}$.\\

Define $f(\beta) = \frac1{n}\sum\limits_{i=1}^{n}(y_{i}e^{-\frac1{2}x_{i}^{T}\beta}+e^{\frac1{2}x_{i}^{T}\beta}) $, from the KKT optimality condition, we know $\lambda\ge \|\nabla f(\beta)\|_\infty$.
Consider the gradient of $f(\beta)$ valued at $\beta=\beta^*$:
\Be\label{eqn:gradient}
\nabla f(\beta^*) = -\frac1{2n}\sum\limits_{i=1}^{n}x_{i}(y_{i}-e^{x_{i}^{T}\beta^*})/\sqrt{e^{x_{i}^{T}\beta^*}},\Ee  and denote $H= \|\nabla f(\beta^*)\|_\infty$.
We will choose a suitable $\lambda$ such that it is greater than $cH$ with high probability and with such a choice of $\lambda$, the estimator is good in the sense that $\|\hat \beta - \beta^*\|_1$ is bounded by a small value which goes to 0 when $n\rightarrow \infty$ under mild conditions.

%Inspired by \cite{belloni2011square-root}, we choose the penalty level $\lambda$ as the following criterion:

%\begin{align}
%{\rm exact\ choice:\ }\lambda&={c}H(1-\alpha|X)\label{eqn:exact}\\
%{\rm asymptotic\ choice:\ }\lambda&=\frac{c}{2}(\sqrt{n})^{-1}\Phi^{-1}(1-\frac{\alpha}{2p})\label{eqn:asymptotic}
%\end{align}
%where $\alpha\in(0,1)$ and $H(1-\alpha|X)$ is the $1-\alpha$ quantile of $H|X$. It is easy to know that $\PP(\lambda\ge cH)\ge1-\alpha$ with $\lambda$ choosing by (\ref{eqn:exact}). While we will prove that under some mild condition the asymptotic choice of $\lambda$ in (\ref{eqn:asymptotic}) yields $\lambda\ge cH$ with probability at least $1-\alpha(1+o(1))$.

In the next section, we study the statistical performance of our proposed method, including how to select a good $\lambda$ such that our estimator has small errors.

\section{Finite-sample and asymptotic bounds on $\|\hat{\beta}-\beta^*\|_1$}
\label{sec: theories}

We first show that when the tuning parameter $\lambda$ defined in our new penalized method \eqref{eqn:mainEQ} is greater than $cH = c\|\nabla f(\beta^*)\|_\infty$ for some $c>1$, the estimation error can be bounded. Based on this result, we will prove a finite sample result for selecting a good tuning parameter  $\lambda$ such that the estimator has very small error with high probability. For our theoretical analysis, we give a few regularity conditions as the following:
\begin{enumerate}[(I)]
\item There exists some positive constant $R<\infty$ such that $\sup\limits_{i\in[n],j\in[p]}|x_{ij}|\le R$.
\item  $n,p$ satisfy that $n\le p \le o(e^{n^{1/5}})$ and $p/\alpha>8$ for $\alpha\in(0,1)$.
\item  For any $p$-dimensional vector $\delta$ in the set $\{\delta\in\R^{p}: \|\delta_{T^{c}}\|_{1}\le L\|\delta_{T}\|_{1}\ {\rm for\ some}\ L>1\}$, there exists some constant $\kappa>0$ such that $\left\langle\delta, \nabla^{2}f(\beta^{*})\delta\right\rangle \ge\kappa^{2}\|\delta_{T}\|_{2}^{2}$.
\end{enumerate}

Conditions (I) and (II) are considerably mild, while (III) is a restricted eigenvalue condition \citep{bickel2008simultaneous}, similar to compatibility condition \citep{Geer2007The} and restricted strong convexity conditions \citep{Negahban2010A}. Although these conditions usually cannot be verified from data, researchers found that they are not strong for Lasso problems in the linear model case and hold with high probability when the elements of design matrix are randomly from Gaussian distributions \citep{Raskutti2010Restricted}. For Poisson regression models, there are a few literatures on when condition (III) holds, we conjecture that it holds with high probability under very mild conditions and leave this to future study. Now we give a deterministic result on when we can have a good estimator in the sense that the  error could be bounded.

\BT \label{thm:deterministic}
Let $\hat{\beta}$ be the estimator defined by \eqref{eqn:mainEQ}. Suppose that the assumptions (I), (II) are satisfied. For some constant $c>1$, assumption (III) is satisfied with $L=\frac{c+1}{c-1}$. If $\lambda$ is chosen such that $\lambda > cH = c\|\nabla f(\beta^*)\|_\infty$ and $\lambda s\le\frac{2\kappa^{2}}{3L(1+L)R}$ then
\begin{align}
\|\hat{\beta}-\beta^{*}\|_{1}&\le\frac{CL(1+L)}{\kappa^{2}}\lambda s,\label{ineq:L1-norm}\\
|f(\hat{\beta})-f(\beta^{*})|&\le\frac{CL(1+L)}{\kappa^{2}}\lambda^{2} s,\label{ineq:error}
\end{align}
for some constant $C\in(2,3]$.
\ET
\BR
Notice that $\|\hat{\beta}-\beta^*\|_2\le\|\hat{\beta}-\beta^*\|_1$. The estimator $\hat{\beta}$ is also $\ell_2$ consistent.
\ER
\vskip 3mm

From Theorem \ref{thm:deterministic}, we see that if we can choose a $\lambda$ such that  $\lambda > cH = c\|\nabla f(\beta^*)\|_\infty$ and $\lambda s\le\frac{2\kappa^{2}}{3L(1+L)R}$  holds with high probability, then conclusions \eqref{ineq:L1-norm} and \eqref{ineq:error} hold with high probability. This motivates us on how to select a good tuning parameter $\lambda$.

Note that $H = \|\nabla f(\beta^*)\|_\infty$  with $f(\beta^*)=-\frac1{2n}\sum\limits_{i=1}^{n}\frac{x_{i}(y_{i}-e^{x_{i}^{T}\beta^*})}{\sqrt{e^{x_{i}^{T}\beta^*}}}$ is a random variable.  Define  by $H(1-\alpha|X)$ the $1-\alpha$ quantile of $H|X$ for  $\alpha\in(0,1)$. If we choose $\lambda$ as follows,
\begin{equation}
\label{eqn:exact}
\lambda=cH(1-\alpha|X)
\end{equation}
 it is easy to know that $\PP(\lambda\ge cH)\ge1-\alpha$ with $\lambda$ choosing by (\ref{eqn:exact}). By a careful analysis, we can prove that $\lambda$ in (\ref{eqn:exact}) is order of $\sqrt{\frac{\log{(p/\alpha)}}{n}}$, which means that $\lambda \leq {\rm const.} \times \sqrt{\frac{\log{(p/\alpha)}}{n}}$ for some ${\rm const.} > 0$.
We shall prove in the appendix the following lemma:

\BL\label{lem:exact}
(i) If $\lambda$ is chosen as (\ref{eqn:exact}), then it implies that $\PP(\lambda\ge cH)\ge1-\alpha$.\\
(ii) Suppose the assumption (I) and (II) are satisfied. Then, $\lambda \leq {\rm const.} \times \sqrt{\frac{\log{(p/\alpha)}}{n}}$ for some $const > 0$.
\EL

\BR
The second conclusion in Lemma \ref{lem:exact}  justifies the condition $\lambda s\le\frac{2\kappa^{2}}{3L(1+L)R}$ if $s\sqrt{\frac{\log p}{n}} \rightarrow 0$.
\ER
\vskip 3mm

Although the $\lambda$ defined in Equation \eqref{eqn:exact} satisfies good property that $\lambda > c H$ and $\lambda s\le\frac{2\kappa^{2}}{3L(1+L)R}$, we notice that this $\lambda$ can not be determined from data in practice because the distribution of $H$ still depends on $\beta^*$. Note that the quantity $\frac{y_{i}-e^{x_{i}^{T}\beta^*}}{\sqrt{e^{x_{i}^{T}\beta^*}}}, i=1,2,\ldots,n$ are i.i.d random variables with mean $0$ and variance $1$, one can approximate choice of $\lambda$ by $\lambda = c\tilde H(1-\alpha|X)$, where $\tilde H(1-\alpha|X)$ is the $1-\alpha$ quantile of $\tilde{H}|X$ and $\tilde{H} = \|-\frac1{2n}\sum\limits_{i=1}^{n}x_{i}z_i\|_\infty$ with $z_i,i=1,2,\ldots, n$ i.i.d.\ from $N(0,1)$. $H$ should have a limiting distribution that is the same as $\tilde H$ under mild conditions.

Motivated by the limiting normal distributions, we can give an asymptotic choice of $\lambda$ such that $\lambda \geq cH$ with high  probability when $n\rightarrow \infty$ as the following:

\begin{equation}
\label{eqn:asymptotic}
\lambda=\frac{c}{2}(\sqrt{n})^{-1}\Phi^{-1}(1-\frac{\alpha}{2p}),
\end{equation}
where $\Phi(\cdot)$ is the cumulative distribution function of the standard norm random variable and $\Phi^{-1}(\cdot)$ is its inverse function. This choice of $\lambda$ has the following properties:

\BL\label{lem:asymp}
(i)  If $\lambda$ is chosen as (\ref{eqn:asymptotic}), \Bes
\begin{array}{rl}
\PP(\lambda\ge cH)&\ge 1-\alpha\left(1+O(1)(\sqrt{2\log{(2p/\alpha)}}-\sqrt{n}b)^{3}n^{-1/2}(3K_{1}\log{p}+b)\right)\\
&\quad\times(1+\frac1{\log{(p/\alpha)}})\frac{\exp\{-2(n\log{(p/\alpha)})^{1/2}b+nb^{2}\}}{1-\sqrt{n}b/(\log{(p/\alpha)})^{1/2}}+C_{1}n/p^{2},
\end{array}
\Ees
where $b=6C_{1}K_{1}\log{p}/p^{3}$ with some positive constants $C_{1}$ and $K_{1}$. In particular, as $n,p\rightarrow\infty$, we have
$$\PP(\lambda\ge cH)\ge 1-\alpha(1+o(1)).$$
(ii) $\lambda$ in (\ref{eqn:asymptotic}) is order of $\sqrt{\frac{\log{(p/\alpha)}}{n}}$.
\EL

Together with Theorem \ref{thm:deterministic} and Lemmas \ref{lem:exact} and \ref{lem:asymp}, we have the following non-asymptotic results.

\BT[Finite-sample]\label{thm:consistency}
Let $\hat{\beta}$ be the estimator defined by \eqref{eqn:mainEQ}. Suppose that assumptions (I), (II) are satisfied. For some constant $c>1$, assumption (III) is satisfied with $L=\frac{c+1}{c-1}$.\\
(i) If $\lambda$ is chosen as (\ref{eqn:exact}) with the above $c$ and the condition $\lambda s\le\frac{2\kappa^{2}}{3L(1+L)R}$ holds, then with probability at least $1-\alpha$, the above inequalities (\ref{ineq:L1-norm}) and (\ref{ineq:error}) hold for some constant $C\in(2,3]$.\\
(ii) If $\lambda$ is chosen as (\ref{eqn:asymptotic}) with the above $c$ and the condition $\lambda s\le\frac{2\kappa^{2}}{3L(1+L)R}$ holds, then for large enough $n$, with probability at least $1-\alpha\left(1+O(1)\frac{\left(\log \frac{p}{\alpha}\right)^{3/2} \log p}{n^{1/2}} \right)$, the above inequalities (\ref{ineq:L1-norm}) and (\ref{ineq:error}) hold for some constant $C\in(2,3]$.
\ET
By Theorem \ref{thm:consistency}, we can obtain the following asymptotic results (or consistency results) immediately as stated in Corollary \ref{coro:asymptotic}.

\BC[Asymptotic]\label{coro:asymptotic}
Let $\hat{\beta}$ be the estimator defined by \eqref{eqn:mainEQ}. Suppose that the assumptions (I), (II) are satisfied. For some constant $c>1$, assumption (III) is satisfied with $L=\frac{c+1}{c-1}$.\\
(i) If $\lambda$ is chosen as (\ref{eqn:exact}) with the above $c$ and the condition $s(\sqrt{n})^{-1}\sqrt{\log{(p/\alpha)}}\rightarrow 0$ holds, then with probability at least $1-\alpha$, the inequalities (\ref{ineq:L1-norm}) and (\ref{ineq:error}) hold for some constant $C\in(2,3]$.\\
(ii) If $\lambda$ is chosen as (\ref{eqn:asymptotic}) with the above $c$ and the condition $s(\sqrt{n})^{-1}\sqrt{\log{(p/\alpha)}}\rightarrow 0$ holds, then with probability at least $1-\alpha(1+o(1))$, the inequalities (\ref{ineq:L1-norm}) and (\ref{ineq:error}) hold for some constant $C\in(2,3]$.
\EC
\BR
The condition $s(\sqrt{n})^{-1}\sqrt{\log{(p/\alpha)}}\rightarrow 0$ means that $s=o(\sqrt{\frac{\log{p}}{n}})$.
\ER

\section{Experiments}
\label{sec:experiments}

We use the R package ``lbfgs" to solve $\ell_1$ penalized convex optimization problems \citep{CoppolaLbfgs}. The lbfgs package implements both the Limited-memory Broyden-Fletcher-Goldfarb-Shanno (L-BFGS) and the Orthant-Wise Quasi-Newton Limited-Memory (OWL-QN) optimization algorithms. The L-BFGS algorithm solves the problem of minimizing an objective, given its gradient, by iteratively computing approximations of the inverse Hessian matrix. The OWL-QN algorithm finds the optimum of an objective plus the L1-norm of the problem's parameters. The package offers a fast and memory-efficient implementation of these optimization routines, which is particularly suited for high-dimensional problems.

\begin{figure}[h!]
\centering
\caption{Success Rates for converge for glmnet and our proposed method (LPWS).  $\beta \sim N(0,1)*c$. We change $c$ from $0.1$ to $3.0$. We did 100 repetitions and the numbers of success of convergence for both algorithms are shown here.  \label{fig:SucRates}}
\includegraphics[scale = 0.6]{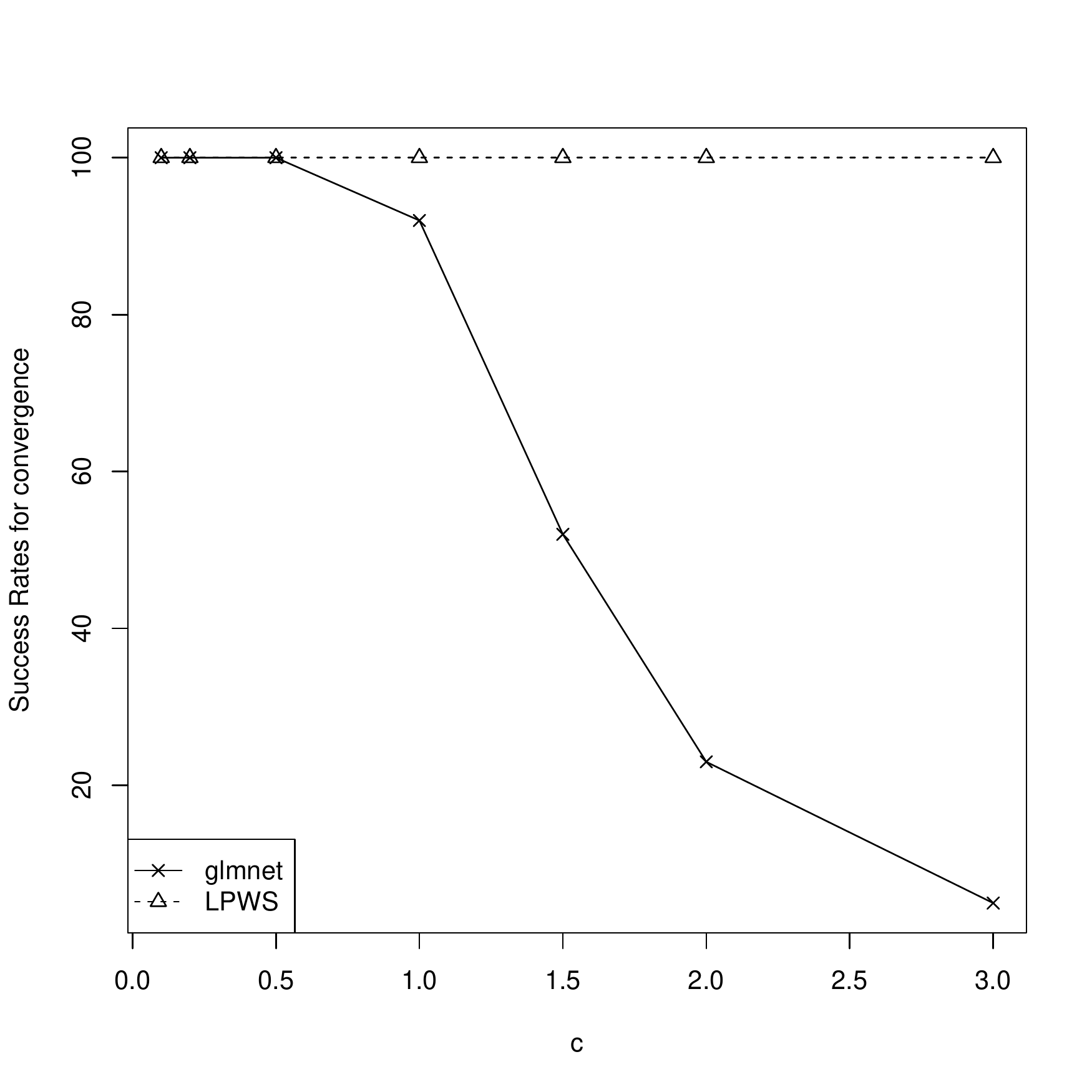}
\end{figure}

We first use simulations to show that our proposed method is much more robust than traditional sparse Poisson models using $\ell_1$ penalized log-likelihood method.  For this purpose, we first generate a design matrix $X \in \mathbb R^{n\times p}$ with $n = 500$, $p = 20$ and each element $X_{ij}$ i.i.d.\ from the standard normal distribution. Then we do centralization and normalization $X$ such that $\sum_{j=1}^n X_{ij} = 0$, and $\frac{1}{n}\sum_{j=1}^n X_{ij}^2= 1, i =1,2,\ldots,n.$ We set the number of nonzero elements of $\beta^*$ as $5$ and each element randomly from $N(0,1)$. $Y_i \sim Poisson(\exp\{\sum_{j=1}^5 X_{ij}\beta_j^*\}).$ We first use R package ``glmnet" to solve the sparse Poisson regression which returns $\ell_1$ regularized log-likelihood estimator. For our proposed method, we set $\lambda=\frac{c}{2}(\sqrt{n})^{-1}\Phi^{-1}(1-\frac{\alpha}{2p})$, with $c = 2$. We repeat this simulation $100$ times and find that there are about 20 times glmnet does not converge and gives warning message or error messages, while our proposed method always converges. If we increase $\beta^*$, glmnet fails more. Below, we provide a plot showing successful convergence rates of glmnet and our proposed method, L1-penalized weighted score (LPWS) method, when the 5 nonzero coefficients are generated from $N(0,1)\times c$ and $c$ varies in the set of $\{0.1,0.2,0.5,1.0,1.5,2.0,3.0\}$. From Figure \ref{fig:SucRates}, we see clearly that our proposed method is much more robust in the sense that it always converges.

\begin{figure}[h!]
\centering
\caption{The solutions of glmnet and our proposed method (LPWS). For glmnet, $\lambda$ is tuned via cross-validation and for LPWS $\lambda=(\sqrt{n})^{-1}\Phi^{-1}(1-\frac{\alpha}{2p})$.  \label{fig:solutions}}
\includegraphics[scale = 0.5]{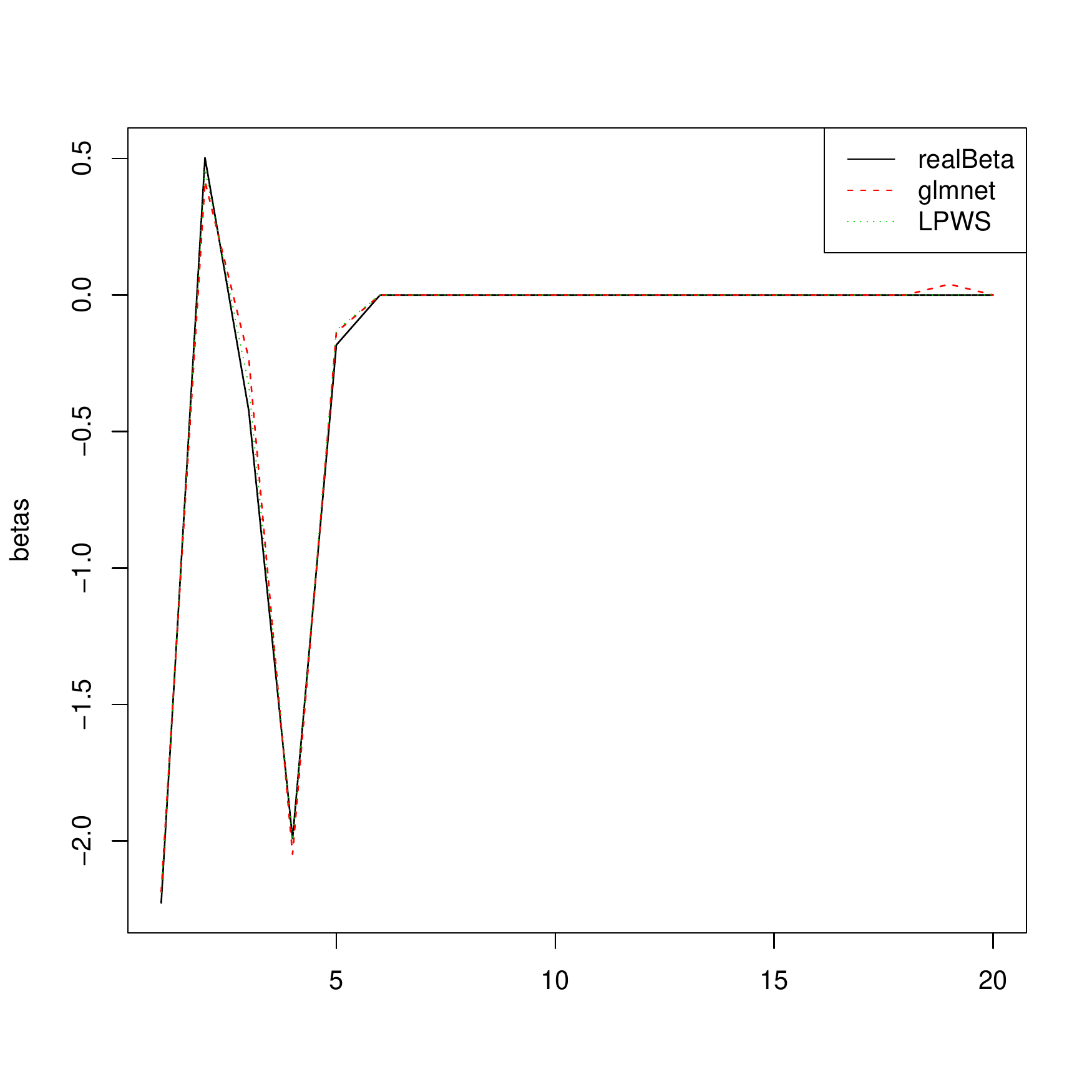}
\end{figure}

To validate the solution of our proposed method using ``lbfgs" package gives good estimator, we use the above simulation settings and choose the nonzero elements of $\beta^*$ from $N(0,1)$. For glmnet, we use cross validation to set the tuning parameter and for our proposed method (LPWS), we choose $\lambda=(\sqrt{n})^{-1}\Phi^{-1}(1-\frac{\alpha}{2p})$.
The solutions and the real coefficients are plotted in Figure \ref{fig:solutions}, from which we see that our new estimator is also a good one. To evaluate the accuracy of our new estimator we do more simulation experiments below.

\begin{figure}[h!]
\centering
\caption{The errors for different methods. ``Err$\_$glmnet" denotes the errors for estimators with glmnet and tuning parameter is selected via cross-validation. ``Err$\_$new" is for our proposed method with $\lambda$ defined as $\lambda=(\sqrt{n})^{-1}\Phi^{-1}(1-\frac{\alpha}{2p})$. ``Err$\_$opt" is for our proposed method  with exact selection of $\lambda$ and ``Err$\_$approx" is for the new proposed method  with an approximate of the exact selection.  \label{fig:errors}}
\includegraphics[scale = 0.5]{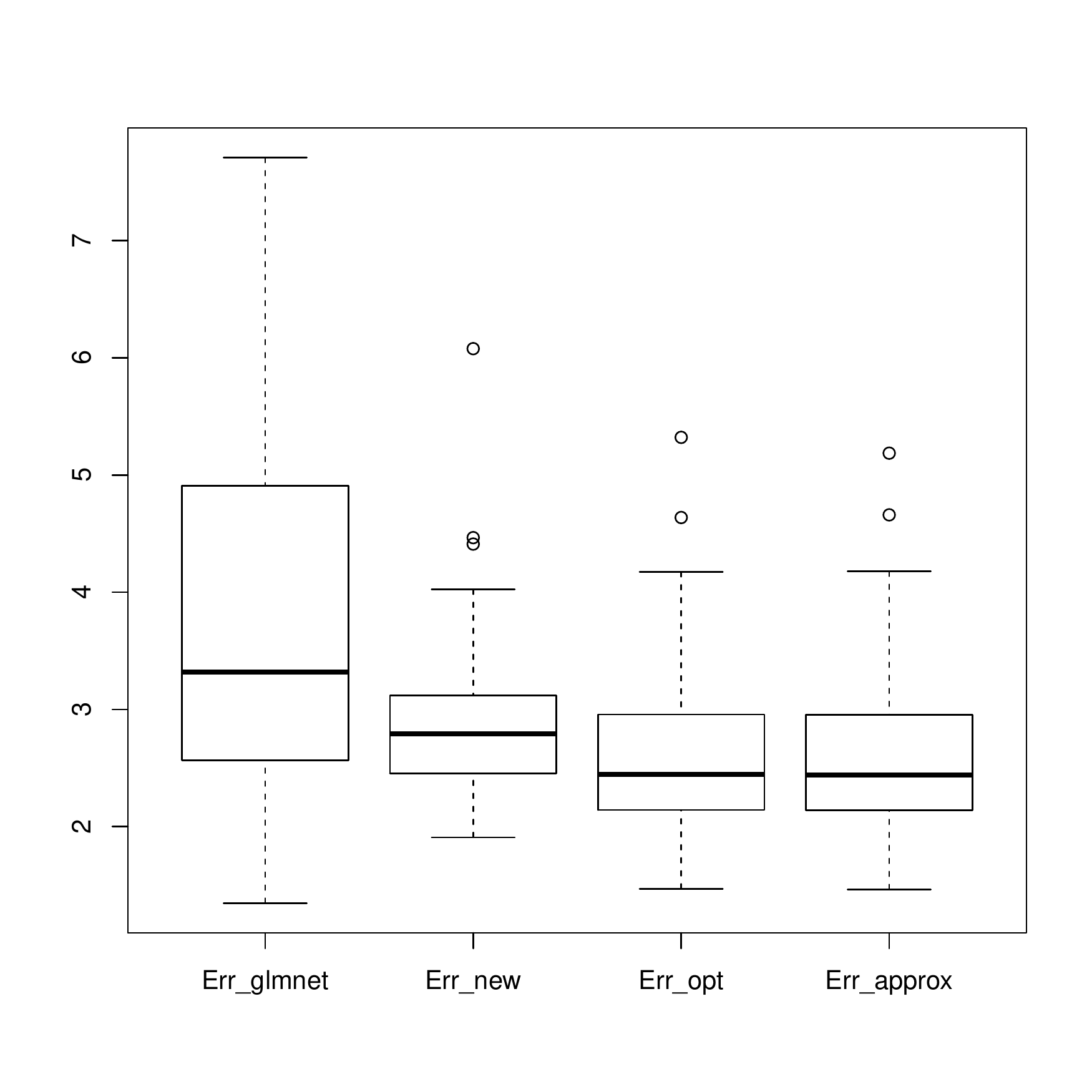}
\end{figure}

Finally, we compare different ways of tuning parameter selection for our proposed method. Recall that $H = \|\nabla f(\beta^*)\|_\infty$. We compare THREE different ways of selecting $\lambda$. (1) As defined in \eqref{eqn:exact}, $\lambda = H(1-\alpha)$. This tuning depends on the real $\beta^*$, which is unknown when we analyze real data, but we still list it here as a benchmark. (2) As defined in \eqref{eqn:asymptotic}, we choose $\lambda = (\sqrt{n})^{-1}\Phi^{-1}(1-\frac{\alpha}{2p})$, this is the asymptotic selection of tuning parameter. (3) We use normal approximation of $H$ defined as $\tilde H = -\frac1{2n}\sum\limits_{i=1}^{n}x_{i}z_i$ with $z_i,i=1,2,\ldots, n$ i.i.d.\ from $N(0,1)$, and define $\lambda = \tilde H(1-\alpha)$. This is an approximation of the exact selection of tuning parameter defined in  \eqref{eqn:exact}.
For comparison, we also calculate the solution of glmnet with $\lambda$ selected using cross validation. In this simulation study, the simulate procedure is almost the same as the previous examples, but here we choose $n = 100, p = 1000$ and  $s = 10$.  We repeat the simulation $100$ times. For each time, we choose each non-zero element of $\beta^*$ from $N(0,1)*0.5$ to make sure that glmnet converges. We calculate the $\ell_1$ estimation error defined as $\|\hat \beta - \beta^*\|_1$. The errors are reported in Figure \ref{fig:errors}, from which we see again that our proposed method outperforms the traditional $\ell_1$ penalized loglikelihood method for sparse Poisson regression,  at the same time, our new method does not need heavy procedure like cross-validation. Hence, our pre-specified tuning parameter works.

\section{Proofs}
\label{sec:proofs}
\noindent
\textbf{Proof of Theorem \ref{thm:deterministic}.}
Let $\delta = \hat \beta - \beta^*$. Recall that $T = \{j:\beta^*_j \neq 0\}$. By definition of $\hat{\beta}$, we have
\Be\label{ineq:bydefinition}
f(\hat{\beta})-f(\beta^{*})\le\lambda(\|\beta^{*}\|_{1}-\|\hat{\beta}\|_{1})
= \lambda [(\|\beta^*_T\|_1 - \|\hat \beta_T\|_1) +(\|\beta^*_{T^c}\| - \|\hat \beta_{T^c}\|_1]
\le\lambda(\|\delta_{T}\|_{1}-\|\delta_{T^{c}}\|_{1}).
\Ee
Since $f(\beta)$ is a convex function, we have
\Be\label{ineq:convexity}
f(\hat{\beta})-f(\beta^{*})\ge\delta^{T}\nabla f(\beta^{*})\ge-\|\nabla f(\beta^{*})\|_{\infty}\|\delta\|_{1}\ge-\frac{\lambda}{c}\|\delta\|_{1},
\Ee
where the last inequality used the choice of $\lambda$ such that $\lambda > cH = c\|\nabla f(\beta^*)\|_\infty$. Combining (\ref{ineq:bydefinition}) and (\ref{ineq:convexity}), we obtain that
\Bes
\lambda(\|\delta_{T}\|_{1}-\|\delta_{T^{c}}\|_{1})\ge-\frac{\lambda}{c}(\|\delta_T\|_{1}+\|\delta_{T^c}\|_{1}),
\Ees
i.e. $$\|\delta_{T^{c}}\|_{1}\le \frac{c+1}{c-1}\|\delta_{T}\|_{1}= L\|\delta_{T}\|_{1}.$$

Defining a new function $\tilde{f}(t)=f(\beta^*+tv)$ from $\R$ to $\R$ for any vector $v\in\R^{p}$, we compute the second and third order derivatives below\\
\Bes
\begin{array}{rl}
\tilde{f}''(t)&=\frac1{2n}\sum\limits_{i=1}^{n}(x_{i}^{T}v)^{2}(y_{i}e^{-x_{i}^{T}(\beta^*+tv)/2}+e^{x_{i}^{T}(\beta+tv)/2})\\
\tilde{f}'''(t)&=-\frac1{4n}\sum\limits_{i=1}^{n}(x_{i}^{T}v)^{3}(y_{i}e^{-x_{i}^{T}(\beta^*+tv)/2}-e^{x_{i}^{T}(\beta^*+tv)/2}).
\end{array}
\Ees
It is easy to obtain that $|\tilde{f}'''(t)|\le\frac1{2}\sup\limits_{i\in[n]}|x_{i}^{T}v|\tilde{f}''(t)\le\frac1{2}\sup\limits_{i\in[n],j\in[p]}|x_{ij}|\|v\|_{1}\tilde{f}''(t)\le \frac1{2}R\|v\|_{1}\tilde{f}''(t)$. Setting $v=\delta=\hat{\beta}-\beta^{*}$, we have
\Be\label{ineq:self-corcondient}
|\tilde{f}'''(t)|\le \frac1{2}R(1+L)\|\delta_{T}\|_{1}\tilde{f}''(t)\le \frac1{2}R(1+L)\sqrt{s}\|\delta_{T}\|_{2}\tilde{f}''(t).
\Ee
Denote $\tilde{R}=\frac1{2}R(1+L)\sqrt{s}$, by Proposition 1 of \cite{Bach2010Self} and (\ref{ineq:self-corcondient}), we have
\Be\label{ineq:taylor}
\begin{split}
f(\hat{\beta})-f(\beta^{*})&\ge\delta^{T}\nabla f(\beta^{*})+\frac{\delta^{T}\nabla^2 f(\beta^{*})\delta}{\tilde{R}^{2}\|\delta_{T}\|_{2}^{2}}(e^{-\tilde{R}\|\delta_{T}\|_{2}}+\tilde{R}\|\delta_{T}\|_{2}-1)\\
&\ge-\|\nabla f(\beta^{*})\|_{\infty}\|\delta\|_{1}+\frac{\delta^{T}\nabla^2 f(\beta^{*})\delta}{\tilde{R}^{2}\|\delta_{T}\|_{2}^{2}}(e^{-\tilde{R}\|\delta_{T}\|_{2}}+\tilde{R}\|\delta_{T}\|_{2}-1)\\
&\ge-\frac{\lambda}{c}\|\delta\|_{1}+\frac{\delta^{T}\nabla^2 f(\beta^{*})\delta}{\tilde{R}^{2}\|\delta_{T}\|_{2}^{2}}(e^{-\tilde{R}\|\delta_{T}\|_{2}}+\tilde{R}\|\delta_{T}\|_{2}-1)
\end{split}
\Ee
Combining (\ref{ineq:bydefinition}) and (\ref{ineq:taylor}), we have
\Be\label{ineq:1}
\begin{split}
&\quad\frac{\delta^{T}\nabla^2 f(\beta^{*})\delta}{\tilde{R}^{2}\|\delta_{T}\|_{2}^{2}}(e^{-\tilde{R}\|\delta_{T}\|_{2}}+\tilde{R}\|\delta_{T}\|_{2}-1)\\
&\le\lambda\|\delta_T\|_1+\frac{\lambda}{c}\|\delta\|_{1}\le L\lambda\|\delta_{T}\|_{1}\le L\lambda\sqrt{s}\|\delta_{T}\|_{2}.
\end{split}
\Ee
By assumption A2 and (\ref{ineq:1}), we have
\Be\label{ineq:2}
e^{-\tilde{R}\|\delta_{T}\|_{2}}+\tilde{R}\|\delta_{T}\|_{2}-1\le\frac{L\lambda \sqrt{s}\tilde{R}^{2}}{\kappa^{2}}\|\delta_{T}\|_{2}.
\Ee
Set
\Be\label{eqn:h}
h=\frac{L\lambda \sqrt{s}\tilde{R}}{\kappa^{2}} = \frac{\lambda RL(1+L)s}{2\kappa^2},
\Ee
then according to the condition on $\lambda$ such that $\lambda s\le\frac{2\kappa^{2}}{3L(1+L)R}$,  we have $h\le\frac1{3}$. Denote $w=\tilde{R}\|\delta_{T}\|_{2}$, then to solve (\ref{ineq:2}) is equivalent to solve the inequality $e^{-w}+w-1\le hw$. By Taylor formula, we have $\frac{w^2}{2}-\frac{w^3}{6}\le e^{-w}+w-1\le hw$ which implies $\{w:e^{-w}+w-1\le hw, h\le \frac1{3}\}\subseteq\{w:\frac{w^2}{2}-\frac{w^3}{6}\le hw, h\le \frac1{3}\}$. Since under the condition $h\le \frac1{3}$, the solution of inequality $\frac{w^2}{2}-\frac{w^3}{6}\le hw$ is $w\le Ch$ for some constant $C\in(2,3]$, then $$\{w:e^{-w}+w-1\le hw, h\le \frac1{3}\}\subseteq\{w\le Ch {\rm\ with\ some\ constant\ }C\in(2,3]\}.$$ So, from (\ref{ineq:2}), we obtain
\Bes
\tilde{R}\|\delta_{T}\|_{2}\le\frac{CL\lambda \sqrt{s}}{\kappa^{2}}\tilde{R},
\Ees
that is,
\Be\label{ineq:4}
\|\delta_{T}\|_{2}\le\frac{CL\lambda \sqrt{s}}{\kappa^{2}}.
\Ee
Hence, notice the relationship $\|\delta\|_{1}\le(1+L)\sqrt{s}\|\delta_{T}\|_{2}$, by (\ref{ineq:4}) we have
\Be\label{ineq:5}
\|\delta\|_{1}\le\frac{CL(1+L)}{\kappa^{2}}\lambda s.
\Ee
Furthermore, by (\ref{ineq:bydefinition}) and (\ref{ineq:convexity}), we obtain
\Be\label{ineq:6}
|f(\hat{\beta})-f(\beta^{*})|\le\lambda\|\delta\|_1\le\frac{CL(1+L)}{\kappa^{2}}\lambda^{2} s.
\Ee

%%%%%%%%%%%%%%%%%%%%%%%%%%%%%%%%%%%%%%%%%%%%%%%%%%%%%%%%%%%%%%%%%%%%%%%%%%%%%%%%%%%%%%%%%%%%%%%%%%%%%%%%

%%%%%%%%%%%%%%%%%%%%%%%%%%%%%%%%%%%%%%%%%%%%%%%%%%%%%%%%%%%%%%%%%%%%%%%%%%%%%%%%%%%

\textbf{Proof of Lemma \ref{lem:exact}}.
(i). By the definition of quantile, it is easy to obtain that
$$\PP(cH>\lambda)=\PP(cH>cH(1-\alpha|X))<\alpha.$$
Then $\PP(\lambda\ge cH)\ge1-\alpha$.\\
(ii). If there exists $t_n=O(\sqrt{\frac{\log{(p/\alpha)}}{n}})$ such that
$\PP(H>t)<\alpha$,
then by definition of quantile we have $H(1-\alpha|X)\le t_n$. So to get $\lambda=cH(1-\alpha|X)=O(\sqrt{\frac{\log{(p/\alpha)}}{n}})$, it suffices to prove that there exists $t_n=O(\sqrt{\frac{\log{(p/\alpha)}}{n}})$ such that
$\PP(H>t_n)<\alpha.$ Let $t=\Phi^{-1}(1-\frac{\alpha}{4p})$ and $t_n=(\sqrt{n})^{-1}t$. It is obvious $t_n=O(\sqrt{\frac{\log{(p/\alpha)}}{n}})$.
Then we shall show that $$\PP(H>t_n)<\frac{\alpha}{2}(1+o(1))<\alpha,$$
as $n,p\rightarrow\infty$  with $n\le p \le o(e^{n^{1/5}})$. \\Recall $H=\max\limits_{j\in[p]}|\frac1{n}\sum\limits_{i=1}^{n}x_{ij}\epsilon_{i}|$ and $\epsilon_{i}=(y_{i}-e^{x_{i}^{T}\beta^*})/\sqrt{e^{x_{i}^{T}\beta^*}}$, then
\Be\label{ineq:prob(i)}
\begin{array}{rl}
\PP(H>t_{n})&=\PP(\max\limits_{j\in[p]}|\frac1{n}\sum\limits_{i=1}^{n}x_{ij}\epsilon_{i}|>(\sqrt{n})^{-1}t)\\
&\le p\max_{j\in[p]}\PP(|\sum\limits_{i=1}^{n}x_{ij}\epsilon_{i}|>\sqrt{n}t).
\end{array}
\Ee
Repeating the argument below (\ref{ineq:prob}), we get

\Bes
\begin{array}{rl}
\PP(H>t_{n})&\le\frac{\alpha}{2}(1+O(1)(\sqrt{2\log{(4p/\alpha)}}-\sqrt{n}b)^{3}n^{-1/2}(3K_{1}\log{p}+b))\\
&\quad\times(1+\frac1{\log{(2p/\alpha)}})\frac{\exp\{-2(n\log{(2p/\alpha)})^{1/2}b+nb^{2}\}}{1-\sqrt{n}b/(\log{(2p/\alpha)})^{1/2}}+C_{1}n/p^{2},
\end{array}
\Ees
where as $n,p\rightarrow\infty$ with $n\le p \le o(e^{n^{1/5}})$, notice that $b$, $\sqrt{n}b$ and $nb^{2}$ are $o(n^{-2})$, we have
$$\PP(H>t_{n})\le\frac{\alpha}{2}(1+o(1))<\alpha.$$
Thus, $\lambda=cH(1-\alpha|X)=O(\sqrt{\frac{\log{(p/\alpha)}}{n}})$.\\

%%%%%%%%%%%%%%%%%%%%%%%%%%%%%%%%%%%%%%%%%%%%%%%%%%%%%%%%%%%%%%%%%%%%%%%%%%%%%%%%%%

\noindent
\textbf{Proof of Lemma \ref{lem:asymp}}.
(i). From (\ref{eqn:gradient}), we know that
$\nabla f(\beta^{*}) = -\frac1{2n}\sum\limits_{i=1}^{n}x_{i}(y_{i}-e^{x_{i}^{T}\beta^*})/\sqrt{e^{x_{i}^{T}\beta^*}}$.
Denote $\epsilon_{i}=(y_{i}-e^{x_{i}^{T}\beta^*})/\sqrt{e^{x_{i}^{T}\beta^*}}$, then
$$H= \|\nabla f(\beta^{*}) \|_{\infty}=\max\limits_{j\in[p]}|\frac1{2n}\sum\limits_{i=1}^{n}x_{ij}\epsilon_{i}|.$$
Denote  $a=\Phi^{-1}(1-\frac{\alpha}{2p})$, then $\lambda=\frac{c}{2}(\sqrt{n})^{-1}a$. Hence
\Be\label{ineq:prob}
\begin{array}{rl}
\PP(cH>\lambda)&=\PP(\max\limits_{j\in[p]}|\frac1{n}\sum\limits_{i=1}^{n}x_{ij}\epsilon_{i}|>(\sqrt{n})^{-1}a)\\
&\le p\max_{j\in[p]}\PP(|\sum\limits_{i=1}^{n}x_{ij}\epsilon_{i}|>\sqrt{n}a).
\end{array}
\Ee
Since $y_{i}|x_{i}\sim Possion(\mu(x_{i}))$ with $\mu_{i}=\mu(x_{i})=e^{x_{i}^{T}\beta^{*}}$, then $\E(e^{\theta\epsilon_{i}})=\exp\{\mu_{i}e^{\theta/\sqrt{\mu_{i}}}-\mu_{i}-\theta/\sqrt{\mu_{i}}\}$ is a positive constant for all $\theta\in\R$. By the exponential Chebyshev's inequality, we have
\Be\label{ineq:sub-exponential}
\PP(|\epsilon_{i}|>M)<e^{-M/K_{1}}\E(e^{\epsilon_{i}/K_{1}})=C_{1}e^{-M/K_{1}}
\Ee
with some constant $C_{1}=\E(e^{\epsilon_{i}/K_{1}})>0$ and $K_{1}>0$.
Denote $\hat{\epsilon}_{i}=\epsilon_{i}1_{\{|\epsilon_{i}|\le M\}}$ and $\check{\epsilon}_{i}=\epsilon_{i}1_{\{|\epsilon_{i}| > M\}}$. Taking $M=3K_{1}\log{p}$, we have
\Bes
\begin{array}{rl}
\PP(|\sum\limits_{i=1}^{n}x_{ij}\epsilon_{i}|>\sqrt{n}a)&= \PP(|\sum\limits_{i=1}^{n}x_{ij}(\hat{\epsilon}_{i}+\check{\epsilon}_{i})|>\sqrt{n}a,\sup\limits_{i\in[n]}|\epsilon_{i}|\le M)\\
&\quad+\PP(|\sum\limits_{i=1}^{n}x_{ij}(\hat{\epsilon}_{i}+\check{\epsilon}_{i})|>\sqrt{n}a,\sup\limits_{i\in[n]}|\epsilon_{i}|> M)\\
&\le\PP(|\sum\limits_{i=1}^{n}x_{ij}\hat{\epsilon}_{i}|>\sqrt{n}a)+\PP(\sup\limits_{i\in[n]}|\epsilon_{i}|> M).
\end{array}
\Ees
Denote $P_{1}=\PP(|\sum\limits_{i=1}^{n}x_{ij}\hat{\epsilon}_{i}|>\sqrt{n}a)$ and $P_{2}=\PP(\sup\limits_{i\in[n]}|\epsilon_{i}|> M)$, then the above inequality can be written as
\Be\label{ineq:truncation}
\PP(|\sum\limits_{i=1}^{n}x_{ij}\epsilon_{i}|>\sqrt{n}a)\le P_{1}+P_{2}.
\Ee
By inequality (\ref{ineq:sub-exponential}) with $M=3K_{1}\log{p}$, we obtain
that
\Be\label{ineq:P2}
P_{2}\le\sum\limits_{i=1}^{n}\PP(|\epsilon_{i}|> M)\le C_{1}ne^{-3\log{p}}=C_{1}n/p^{3}.
\Ee
To estimate the $P_{1}$, we need the following Sakhanenko type moderate deviation theorem of \cite{Sakhanenko1991Berry}, i.e.
\BL\label{lem:MDT}
Let $\eta_{1},\cdots,\eta_{n}$ be independent random variables with $\E\eta_{i}=0$ and $|\eta_{i}|<1$ for all $i\in[n]$. Denote $\sigma_{n}^{2}=\sum\limits_{i=1}^{n}\E\eta_{i}^{2}$ and $L_{n}=\sum\limits_{i=1}^{n}\E|\eta_{i}|^{3}/\sigma_{n}^{3}$. Then there exists a positive constant $A$ such that for all $x\in[1,\frac1{A}\min\{\sigma_{n},L_{n}^{-1/3}\}]$
\Bes
\PP(\sum\limits_{i=1}^{n}\eta_{i}>x\sigma_{n})=(1+O(1)x^{3}L_{n})\bar\Phi(x),
\Ees
where $\bar\Phi(x)=1-\Phi(x)$ and $\Phi(x)$ is the cumulative distribution function of standard normal distribution.
\EL
Since $\E(\epsilon_{i})=\E(\hat{\epsilon}_{i})+\E(\check{\epsilon}_{i})=0$, then it is easy to obtain that $$|\E(\hat{\epsilon}_{i})|=|\E(\check{\epsilon}_{i})|\le 2C_{1}Me^{-M/K_{1}}=6C_{1}K_{1}\log{p}/p^{3}.$$ Denote $b=6C_{1}K_{1}\log{p}/p^{3}$, then $|\E(\hat{\epsilon}_{i})|\le b$ and $$|\sum\limits_{i=1}^{n}x_{ij}\E\hat{\epsilon}_{i}|\le\sqrt{(\sum\limits_{i=1}^{n}x_{ij}^{2})(\sum\limits_{i=1}^{n}|\E\hat{\epsilon}_{i}|^{2})}\le nb.$$
Furthermore, with assumption A1 we have $$|x_{ij}(\hat{\epsilon}_{i}-\E\hat{\epsilon}_{i})|\le(\sup\limits_{i\in[n],j\in[p]}|x_{ij}|)(|\epsilon_{i}|+|\E\epsilon_{i}|)\le R(M+b).$$ Notice the inequality $\E\hat{\epsilon}_{i}^{2}=\E\epsilon_{i}^{2}-\E\check{\epsilon}_{i}^{2}\le\E\epsilon_{i}^{2}=1$. Let $\eta_{ij}=x_{ij}(\hat{\epsilon}_{i}-\E\hat{\epsilon}_{i})/R(M+b)$, then we have $\E\eta_{ij}=0, |\eta_{ij}|<1$,
\Bes
\begin{array}{rl}
\sigma_{nj}^{2}&=\sum\limits_{i=1}^{n}\E\eta_{ij}^{2}=\frac1{R^{2}(M+b)^{2}}\sum\limits_{i=1}^{n}\E(x_{ij}^{2}(\hat{\epsilon}_{i}-\E\hat{\epsilon}_{i})^{2})\\
&\le\frac1{R^{2}(M+b)^{2}}\sum\limits_{i=1}^{n}x^{2}_{ij}\E\hat{\epsilon}_{i}^{2}\le\frac1{R^{2}(M+b)^{2}}\sum\limits_{i=1}^{n}x_{ij}^{2}\\&=\frac{n}{R^{2}(M+b)^{2}},\\
L_{nj}&=\sum\limits_{i=1}^{n}\E|\eta_{ij}|^{3}/\sigma_{nj}^{3}\le \sum\limits_{i=1}^{n}\E|\eta_{ij}|^{2}/\sigma_{nj}^{3}=\frac1{\sigma_{nj}}.
\end{array}
\Ees
Then, $\sigma_{nj}^{2}=O(\frac{n}{(M+b)^2})$ and $L_{nj}=O(\frac{M+b}{\sqrt{n}})$. Using Lemma \ref{lem:MDT}, for large enough $n,p$ such that $n\le p \le o(e^{n^{1/5}})$, we have
\Be\label{P1}
\begin{array}{rl}
P_{1}&=\PP(|\sum\limits_{i=1}^{n}x_{ij}(\hat{\epsilon}_{i}-\E\hat{\epsilon}_{i}+\E\hat{\epsilon}_{i})|>\sqrt{n}a)\\
&\le\PP(|\sum\limits_{i=1}^{n}x_{ij}(\hat{\epsilon}_{i}-\E\hat{\epsilon}_{i}))|>\sqrt{n}a-|\sum\limits_{i=1}^{n}x_{ij}\E\hat{\epsilon}_{i}|)\\
&\le \PP(|\sum\limits_{i=1}^{n}\frac{x_{ij}(\hat{\epsilon}_{i}-\E\hat{\epsilon}_{i}))}{R(M+b)}|>\frac{\sqrt{n}}{R(M+b)}(a-\sqrt{n}b))\\
&\le \PP(|\sum\limits_{i=1}^{n}\eta_{ij}|>\sigma_{nj}(a-\sqrt{n}b))\\
&=2(1+O(1))(a-\sqrt{n}b)^{3}L_{nj}\bar{\Phi}(a-\sqrt{n}b)\\
\end{array}
\Ee
with $a-\sqrt{n}b$ uniformly in $[1, O(n^{1/6}(\log{p})^{-1/3})]$. Notice that $\log{(p/\alpha)}<a^{2}<2\log{(2p/\alpha)}$ when $p/\alpha>8$ and for all $u>0$ the inequality $\frac{u}{1+u^{2}}\phi(u)\le\bar{\Phi}(u)\le\frac{\phi(u)}{u}$ holds where $\phi(\cdot)$ is the density function of standard normal distribution. Then,
\Be\label{ineq:Phibar}
\begin{array}{rl}
\bar{\Phi}(a-\sqrt{n}b)&\le\frac{\phi(a-\sqrt{n}b)}{a-\sqrt{n}b}=\phi(a)\frac{\exp\{-2a\sqrt{n}b+nb^{2}\}}{a-\sqrt{n}b}\\
&=\frac{a}{1+a^{2}}\phi(a)\frac{1+a^{2}}{a(a-\sqrt{n}b)}\exp\{-2a\sqrt{n}b+nb^{2}\}\\
&\le\bar{\Phi}(a)\frac{1+a^{2}}{a(a-\sqrt{n}b)}\exp\{-2a\sqrt{n}b+nb^{2}\}\\
&= \frac{\alpha}{2p}(1+\frac1{a^{2}})\frac1{1-\sqrt{n}b/a}\exp\{-2a\sqrt{n}b+nb^{2}\}\\
&\le\frac{\alpha}{2p}(1+\frac1{\log{(p/\alpha)}})\frac{\exp\{-2(n\log{(p/\alpha)})^{1/2}b+nb^{2}\}}{1-\sqrt{n}b/(\log{(p/\alpha)})^{1/2}}
\end{array}
\Ee
and
\Be\label{ineq:o(1)}
\begin{array}{rl}
(a-\sqrt{n}b)^{3}L_{nj}\le(\sqrt{2\log{(2p/\alpha)}}-\sqrt{n}b)^{3}n^{-1/2}(3K_{1}\log{p}+b).
\end{array}
\Ee
Combining (\ref{P1}), (\ref{ineq:Phibar}) and (\ref{ineq:o(1)}), we have
\Be\label{ineq:P1}
\begin{array}{rl}
P_{1}&\le\frac{\alpha}{p}(1+O(1)(\sqrt{2\log{(2p/\alpha)}}-\sqrt{n}b)^{3}n^{-1/2}(3K_{1}\log{p}+b))\\
&\quad\times(1+\frac1{\log{(p/\alpha)}})\frac1{1-\sqrt{n}b/(\log{(p/\alpha)})^{1/2}}\exp\{-2(n\log{(p/\alpha)})^{1/2}b+nb^{2}\}.
\end{array}
\Ee
Thus, combining (\ref{ineq:prob}), (\ref{ineq:P2}) and (\ref{ineq:P1}), we obtain that
\Bes
\begin{array}{rl}
&\quad\PP(c\|\nabla f(\beta^{*})\|_{\infty}>\lambda)\le p(P_{1}+P_{2})\\
&\le\alpha(1+O(1)(\sqrt{2\log{(2p/\alpha)}}-\sqrt{n}b)^{3}n^{-1/2}(3K_{1}\log{p}+b))\\
&\quad\times(1+\frac1{\log{(p/\alpha)}})\frac{\exp\{-2(n\log{(p/\alpha)})^{1/2}b+nb^{2}\}}{1-\sqrt{n}b/(\log{(p/\alpha)})^{1/2}}+C_{1}n/p^{2}.
\end{array}
\Ees
As $n,p\rightarrow\infty$ with $n\le p \le o(e^{n^{1/5}})$, notice that $b$, $\sqrt{n}b$ and $nb^{2}$ are $o(n^{-2})$, hence, we have
$$\PP(c\|\nabla f(\beta^{*})\|_{\infty}>\lambda)\le\alpha(1+o(1)).$$
(ii). Notice the fact that for any $u>0$, the inequality $$1-\Phi(u)\le\frac{\phi(u)}{u}$$ holds where the $\phi(\cdot)$ is the density function of standard normal distribution.  Let $u=\Phi^{-1}(1-\frac{\alpha}{2p})$. If $p/\alpha>8$, it is easy to see $u>3/2$. Then the above inequality becomes
\Bes
\frac{\alpha}{2p}=1-\Phi(u)\le\frac{\phi(u)}{u}=\frac{\exp\{-u^{2}/2\}}{\sqrt{2\pi}u}<\exp\{-u^{2}/2\},
\Ees
i.e. $u<\sqrt{2\log{(2p/\alpha)}}$. Thus $\Phi^{-1}(1-\frac{\alpha}{2p})=O(\sqrt{\log(p/\alpha)})$ and
\Bes
\lambda=c(\sqrt{n})^{-1}\Phi^{-1}(1-\frac{\alpha}{2p})=O(\sqrt{\frac{\log{(p/\alpha)}}{n}}).
\Ees
%%%%%%%%%%%%%%%%%%%%%%%%%%%%%%%%%%%%%%%%%%%%%%%%%%%%%%%%%%%%%%%%%%%%%%%%%%%%%%%%%%%%%%

%{\bf Acknowledgements}: FX is supported by the grants Macao S.A.R FDCT  049/2014/A1,  MYRG2015-00021-FST and MYRG2016-00025-FST. LX is supported by the grants NNSFC 11571390, Macao S.A.R FDCT  049/2014/A1, MYRG2015-00021-FST and MYRG2016-00025-FST.

%\begin{thebibliography}{9}
%\bibitem{Bach10} Bach, Francis. Self-concordant analysis for logistic regression. Electron. J. Statist. 4 (2010), 384--414.
%\bibitem{BCW11} Belloni, A., Chernozhukov, V. and Wang, L., Pivotal recoveryof sparse signals via conic programming. Biometrika, Vol. 98, 4(2011),791-806.
%\bibitem{Sak92}Sakhanenko, A. I., Berry-Esseen type estimates for large deviation probabilities.Siberian Math. J. 32, 647-656(1992).
%\bibitem{LV15} Li, Yen-Huan and Volkan Cevher. Consistency of ??-regularized maximum-likelihood for compressive Poisson regression. ICASSP (2015).
%\end{thebibliography}

\end{document}